\documentclass[reqno,12pt]{amsart}

\usepackage{amscd,amssymb}

\theoremstyle{plain}
\newtheorem{Thm}[subsection]{Theorem}
\newtheorem{Cor}[subsection]{Corollary}
\newtheorem{Lem}[subsection]{Lemma}
\newtheorem{Prop}[subsection]{Proposition}

\theoremstyle{definition}
\newtheorem{Def}[subsection]{Definition}

\theoremstyle{remark}

\newtheorem{Rem}[subsection]{Remark}

\numberwithin{equation}{section}
\renewcommand{\rm}{\normalshape}

\newif\ifShowLabels
\ShowLabelstrue
\newdimen\theight
\def\TeXref#1{%
    \leavevmode\vadjust{\setbox0=\hbox{{\tt
        \quad\quad  {\small \rm #1}}}%
    \theight=\ht0
    \advance\theight by \lineskip
    \kern -\theight \vbox to
    \theight{\rightline{\rlap{\box0}}%
    \vss}%
    }}%

\ShowLabelsfalse% comment this out if labels should be printed

\renewcommand{\sec}[2]{\section{#2}\label{S:#1}%
    \ifShowLabels \TeXref{{S:#1}} \fi}
\newcommand{\ssec}[2]{\subsection{#2}\label{SS:#1}%
    \ifShowLabels \TeXref{{SS:#1}} \fi}

\newcommand{\refs}[1]{Section ~\ref{S:#1}}
\newcommand{\refss}[1]{Section ~\ref{SS:#1}}

\newcommand{\reft}[1]{Theorem ~\ref{T:#1}}
\newcommand{\refl}[1]{Lemma ~\ref{L:#1}}
\newcommand{\refp}[1]{Proposition ~\ref{P:#1}}
\newcommand{\refc}[1]{Corollary ~\ref{C:#1}}

\newcommand{\refe}[1]{\eqref{E:#1}}

\newenvironment{thm}[1]%
    { \begin{Thm} \label{T:#1}  \ifShowLabels \TeXref{T:#1} \fi }%
    { \end{Thm} }

\renewcommand{\th}[1]{\begin{thm}{#1} \sl }
\renewcommand{\eth}{\end{thm} }

\newenvironment{lemma}[1]%
    { \begin{Lem} \label{L:#1}  \ifShowLabels \TeXref{L:#1} \fi }%
    { \end{Lem} }
\newcommand{\lem}[1]{\begin{lemma}{#1} \sl}
\newcommand{\elem}{\end{lemma}}

\newenvironment{propos}[1]%
    { \begin{Prop} \label{P:#1}  \ifShowLabels \TeXref{P:#1} \fi }%
    { \end{Prop} }
\newcommand{\prop}[1]{\begin{propos}{#1}\sl }
\newcommand{\eprop}{\end{propos}}

\newenvironment{corol}[1]%
    { \begin{Cor} \label{C:#1}  \ifShowLabels \TeXref{C:#1} \fi }%
    { \end{Cor} }
\newcommand{\cor}[1]{\begin{corol}{#1} \sl }
\newcommand{\ecor}{\end{corol}}

\newenvironment{defeni}[1]%
    { \begin{Def} \label{D:#1}  \ifShowLabels \TeXref{D:#1} \fi }%
    { \end{Def} }
\newcommand{\defe}[1]{\begin{defeni}{#1} \sl }
\newcommand{\edefe}{\end{defeni}}

\newenvironment{remark}[1]%
    { \begin{Rem} \label{R:#1}  \ifShowLabels \TeXref{R:#1} \fi }%
    { \end{Rem} }
\newcommand{\rem}[1]{\begin{remark}{#1}}
\newcommand{\erem}{\end{remark}}

\newcommand{\eq}[1]%
    { \ifShowLabels \TeXref{E:#1} \fi
       \begin{equation} \label{E:#1} }
\newcommand{\eeq}{ \end{equation} }

\newcommand{\prf}{ \begin{proof} }
\newcommand{\epr}{ \end{proof} }

\newcommand\alp{\alpha}     

     \newcommand\Gam{\Gamma}

\newcommand\lam{\lambda}        \newcommand\Lam{\Lambda}

\newcommand\ome{\omega}     

\newcommand\calA{{\mathcal{A}}}
\newcommand\calB{{\mathcal{B}}}
\newcommand\calC{{\mathcal{C}}}
\newcommand\calD{{\mathcal{D}}}
\newcommand\calE{{\mathcal{E}}}

\newcommand\calH{{\mathcal{H}}}

\newcommand\calM{{\mathcal{M}}}

\newcommand\calO{{\mathcal{O}}}
\newcommand\calP{{\mathcal{P}}}

\newcommand\calS{{\mathcal{S}}}

        \newcommand\bfF{{\mathbf F}}

\newcommand\RR{\mathbb{R}}

\renewcommand\AA{\mathbb{A}}

\newcommand\FF{\mathbb{F}}

\newcommand\CC{\mathbb{C}}

 \newcommand\grh{{\mathfrak{h}}}

\newcommand\sdp{\times \hskip -0.3em {\raise 0.3ex
\hbox{$\scriptscriptstyle |$}}} % semidirect product

\newcommand\End{\operatorname{End\,}}
\newcommand\Ext{\operatorname{Ext}}

\newcommand\Hom{\operatorname {Hom}}

\newcommand\Ker{\operatorname{Ker}}
\newcommand\Coker{\operatorname{Coker}}

\newcommand\Perv{\operatorname{Perv}}

\newcommand\Spec{\operatorname{Spec}}

\newcommand\oX{{\overline{X}}}

\newcommand\opi{{\overline{\pi}}}

\newcommand\tilf{{\widetilde{f}}}

\newcommand\tilX{{\widetilde{X}}}

\newcommand\tilY{{\widetilde{Y}}}

\newcommand\Phitil{{\tilde{\Phi}}}

\newcommand\ten{\otimes}

\newcommand{\To}{\longrightarrow}
\newcommand{\iso}{\widetilde\longrightarrow}
\newcommand{\imbed}{\hookrightarrow}
\newcommand{\oplusl}{\bigoplus\limits}
\newcommand{\cupl}{\bigcup\limits}

\newcommand{\can}{ \tilde F^*}

\newcommand{\A}{{\mathcal A}}
\newcommand{\B}{{\mathcal B}}

\newcommand\Id{\text{Id}}

\newcommand\dg{\calD}

\newcommand\M{\calM}
\newcommand\C{\CC}
\newcommand\h{\grh}
\newcommand\galochka{^\vee{\;}}
\newcommand\order{\operatorname{order}}
\begin{document}
\title[]{Gluing of abelian categories and differential operators on
the basic affine space}
\author{Roman Bezrukavnikov, Alexander Braverman, Leonid Positselskii}
\address{Department of Mathematics, University of Chicago, Chicago IL 60637}
\address{Department of Mathematics, Harvard University, 1 Oxford St.,
Cambridge MA, 02138}
\address{Independent Moscow University and IHES}
\email{braval@math.harvard.edu, roman@math.uchicago.edu, posic@ihes.fr}
\thanks{A.B and R.B. are partially supported by the National Science
Foundation}
\begin{abstract}
The notion of gluing of abelian categories was introduced in \cite{kl}
and studied in \cite{Po}. We observe that this notion
is a particular case of a general categorical construction.

We then apply this general notion to the study of the ring of
global differential operators $\calD$ on the  basic affine space $G/U$
(here $G$ is a semi-simple simply connected algebraic group over $\CC$
and $U\subset G$ is a maximal unipotent subgroup).

We show that the category of $\calD$-modules is glued from $|W|$ copies
of the category of $D$-modules on $G/U$ where $W$ is the Weyl group, and
the Fourier transform is used to define the gluing data.
As an application we prove that the algebra $\calD$ is noetherian,
and get some information on its  homological properties.
\end{abstract}
\maketitle
%--------------------------------------------------------
\sec{int}{Introduction}
\ssec{}{}Let $Y$ be an affine algebraic  variety over $\CC$.
Let $\calO$ be the ring of regular functions on $Y$, and  $\calD$ be the
algebra of differential operators on $\calO$ in the sense of Grothendieck
 (thus $\calD=\cupl_n
\calD_n$, where for $d:\calO\to \calO$ we have $d\in \calD_n\iff
[f,d]\in \calD_{n-1}$ for all $f\in \calO$; and $D_{-1}=0$). If $Y$ is
normal, then $\calD$ is identified with  global sections of the
sheaf of differential operators  with polynomial coefficients on
the smooth part $Y^{sm}$.

If $Y=Y^{sm}$ is smooth, then it is well-known that  algebra
$\calD$ possesses the following properties:

a) $\calD$ is noetherian

b) The homological dimension of $\calD$ is equal to $\dim Y$.

On the other hand, it is known after the work of J.~Bernstein, I.~Gelfand
and S.~Gelfand (cf. \cite{bgg2})
 that a) and b) fail if we do not assume that $Y$ is
non-singular. In fact in \cite{bgg3} the authors show that the algebra
$\calD$ does not "behave nicely" already when $Y$ is the cubic cone in a
three-dimensional space.
%-----------------------------------------------------------
\ssec{}{The basic affine space}In this paper we exhibit
some examples of singular affine varieties $Y$ for which the algebra $\calD$
of global differential operators behaves somewhat similarly to the algebra of
differential operators on a non-singular variety. These examples come
from semi-simple groups.

Namely, let $G$ be a semi-simple, simply connected algebraic group
over $\CC$ and let $U\subset G$ be a maximal unipotent subgroup of
$G$. Consider the variety $X=G/U$. Then it is known that $X$ is a
quasi-affine variety. This means that the algebra $\calO$ of
global regular functions on $X$ is finitely-generated and
separates points of $X$. Let $\oX=\Spec \calO$. Then it is easy to
see that $\oX$ is singular unless $G=SL(2)^n$. Let $\calD$ be the
algebra of global differential operators on $X$. %By the definition
Since $\oX$ is normal,
$\calD$ is equal to the algebra of  differential operators
on $\oX$. The following result is proven in \refs{diff}. \th{one}
The algebra $\calD$ is (left and right) noetherian. \eth
\reft{one} is a generalization of a) for $Y=\oX$. Let us now turn
to the generalization of b). The homological dimension of the
algebra $\calD$ is probably infinite, cf. \cite{PB}.
However in \refs{diff} we prove the following
\th{two}The  %homological
injective dimension of $\calD$ as a left $\calD$-module is equal
to $\dim X$. \eth Note that  $\calD$ is a projective generator of
the category of left $\calD$-modules; and for a category of {\em
finite} homological dimension the homological dimension of the category
equals the injective
dimension of a projective generator (thus if the homological
dimension of this category were finite then it would be equal to
$\dim X$). Also \reft{two} implies that one can define an analogue
of Verdier duality for $\calD$-modules.
%-----------------------------------------------------------------------------
%--------
\ssec{}{Gluing}Denote by $D_X$ the sheaf of differential operators on $X$.
Let us now explain the relation between the category
$\calM(\calD)$ of left $\calD$-modules and the category
$\calM(D_X)$ of $D$-modules on $X$ (i.e. quasi-coherent sheaves of
modules over the sheaf of algebras $D_X$).

In \cite{kl} D.~Kazhdan and G.~Laumon considered the category $\Perv(X)$
of $\ell$-adic
perverse sheaves on the variety $X$ over a finite field $\FF_q$. In
particular they introduced for every element $w$ of the Weyl group $W$ of
$G$ a functor $\bfF_w:\Perv(X)\to\Perv(X)$. The functors $\bfF_w$ are
generalizations of the Fourier-Deligne transform. Using the functors
$\bfF_w$ they defined a new {\it glued} category $\calP$. For example,
when $G=SL(2)$ then $X=\AA^2\backslash\{0\}$ and the Kazhdan-Laumon category
 $\calP$
is equivalent to the category of perverse sheaves on $\AA^2$. However, for
general $G$ the category $\calP$ is not equivalent to the category of
perverse sheaves on any variety.
This construction was studied by A.~Polishchuk in \cite{Po}.
%(our definition of gluing is partly borrowed from {\it loc. cit.})

The main result of this paper asserts that the $D$-module
counterpart of the Kazhdan-Laumon construction produces a category
equivalent to $\calM(\calD)$; thus it says that  $\calM(\calD)$
is glued from $|W|$ copies of $\calM(D_X)$.
This result implies both \reft{one} and \reft{two}.

The procedure of Kazhdan-Laumon gluing is in fact a particular
case of a very general categorical construction called the
category of coalgebras over a comonad   (alternatively, it can also
be obtained as a category of algebras over a monad).
%, as existence of both left and right adjoints  is assumed in \cite{kl}).
 For a pair of adjoint functors between two
categories $\A$, $\B$ satisfying some additional assumptions it
allows to describe the category $\A$ in terms of objects in $\B$
endowed with an additional structure. These assumptions hold in
particular if the categories are abelian, and one of the two
functors is exact and faithful.

The same notion of gluing was used by A.~Rosenberg in his paper about
noncommutative schemes~\cite{Ros}.
Kontsevich and Rosenberg~\cite{KR} used comonads of a more general nature
for construction of noncommutative stacks.

The proof of our main result reduces to the following. The algebra
$\calD$ carries an action of the Weyl group $W$, with the action
of a simple reflection defined by means of Fourier transform (this
statement is in fact due to Gelfand and Graev, \cite{GG}). Let
$L:\calM(\calD)\to\calM(D)$ denote the localization functor
sending every $M\in\calM(\calD)$ to $\calO_X\underset{\calO}\ten
M$. Then our gluing theorem is equivalent to the fact that for
every non-zero $M\in\calM(\calD)$ there exists $w\in W$ such that
$L(M^w)\neq 0$. This fact is proved by a direct calculation,
similar to some standard computations in the theory of semi-simple
Lie algebras (e.g. the computation of the determinant of the
Shapovalov form on a Verma module).
%-----------------------------------------------------------------------------
\ssec{ack}{Acknowledgements}This paper ows its very existence to D.~Kazhdan
who communicated to us its main results as a conjecture. The second author
is also grateful to J.~Bernstein and A.~Polishchuk for very useful discussions
on the subject. The third author is grateful to D.~Kaledin, who first
told him about monads many years ago.

The proof of \reft{mainr} was obtained by the first two authors in~1995,
when they were graduate students at the University of Tel-Aviv.
They take the pleasure of acknowledging their gratitude to alma mater.
Likewise, the third author did his part of this research in~1996, when
he was a graduate student at Harvard University.  He is glad to use
the opportunity to express his gratitude to this institution.
%
%------------------------------------------------------------------------------
%----------------------
\sec{abs}{Gluing of abelian categories}

\ssec{mona}{Comonads and coalgebras over comonads}
 In this section we recall some general categorical constructions,
most of which can be found e.g.\ in chapter VI of the book~\cite{ML}
(one has to reverse the direction of various arrows,
i.e., to replace all categories by the opposite ones,
to pass from our setting to the one of {\it loc.\ cit.})
%We sketch some of the arguments for the reader's convenience.
Examples are provided in the end of this section and in \refs{diff}.

Let $\calB$ be a category.  The category of functors from $\calB$
to itself is a monoidal category with respect to the operation of
composition.  Therefore, one can consider monoid or comonoid objects
in this category.  In~\cite{ML} they are called ``monads'' and
``comonads'' in $\B$. Let us spell out the resulting definition.

\defe{gldef}  %(see \cite{ML}, chapter VI)
A comonad in $\calB$ is a functor $\Phi:\calB\to\calB$ together
with morphisms of functors $\eta: \Phi\to \Id$ (counit),
$\mu:\Phi\to \Phi\circ \Phi$ (comultiplication)
such that
the two compositions
\begin{equation}\label{eqmona1}
\begin{array}{ll}
\Phi \overset{\mu}{\To}\Phi^2 \overset{\Phi\mu}{\To}{\Phi^3}\\
\Phi \overset{\mu}{\To}\Phi^2 \overset{\mu\Phi}{\To}{\Phi^3}
\end{array}
\end{equation}
coincide; and each of the two compositions
\begin{equation}\label{eqmona2}
\begin{array}{ll}
\Phi \overset{\mu}{\To}\Phi^2\overset{\Phi\eta}{\To} \Phi\circ \Id=\Phi \\
\Phi \overset{\mu}{\To}\Phi^2\overset{\eta\Phi}{\To} \Id\circ
\Phi=\Phi
\end{array}
\end{equation}
equals identity.
\edefe

In the next definition we stick to the terminology from~\cite{ML}.
In our situation (when $\B$ is abelian and $\Phi$ is additive)
it would be more natural to call such objects ``comodules'' rather
then ``coalgebras''.

\defe{}
A coalgebra over a comonad $(\Phi,\eta,\mu)$
is an object $X$ of $\calB$ together with a morphism $h:X\to \Phi(X)$
(comultiplication)
such that
\begin{equation}\label{eqcoa}
\begin{array}{ll}
\eta_X\circ h=id_X;\\
\mu_X\circ h=\Phi(h) \circ h.
\end{array}
\end{equation}
\edefe

Coalgebras over a comonad $(\Phi, \eta,\mu)$ form a category, which
we denote by $\calB_\Phi$.

%%If $B$ is abelian and $\Phi$ is additive,
%We will also say ``a $\Phi$-module'' instead of ``a coalgebra over
%$(\Phi,\eta,\mu)$'', especially if $\B$ is abelian and $\Phi$ is
%additive.

\lem{24}
Let $\B$ be an abelian category, and $(\Phi, \eta,\mu)$ be a comonad
on $\calB$. Assume that $\Phi$ is (additive and) left exact; then
$\calB_\Phi$ is an abelian category, and the forgetful functor
$\B_\Phi \to \B$ is exact.
\elem

\prf For $\tilX=(X,h_X), \tilY=(Y,h_Y)\in \B_\Phi$ and a morphism
$\tilf:\tilX\to\tilY$ we need to check
that kernel and cokernel of $\tilf$ exist, and that
\begin{equation}\label{proverab}
\Coker(\Ker(\tilf)\to X)\iso
\Ker (Y\to \Coker (\tilf))
\end{equation}
(the other axioms are clear). Let $f$ denote $\tilf$ considered as an
element of $\Hom_\B(X,Y)$. Then there are unique
morphism $h_{\Ker(f)}:\Ker(f)\to \Phi(\Ker(f))$ and $h_{\Coker(f)}:\Coker(f)
\to \Phi(\Coker(f))$ which make the obvious diagrams commutative;
moreover, $(\Ker(f), h_{\Ker(f)})$, $(\Coker(f), h_{\Coker(f)})$ are easily
seen to be a kernel, and a cokernel of $\tilf$
(here the
diagrams involving $\Ker(f)$ are commutative because they inject into
the corresponding diagrams for $X$; while in the ones involving $\Coker(f)$
one only needs to verify equalities of elements in $\Hom(\Coker(f)),Z)$
for various objects $Z\in \B$, thus it is enough to verify the equality
of their compositions with a surjective arrow $Y\to \Coker(f)$).

This shows existence of $\Ker(\tilf)$, $\Coker(\tilf)$; since an arrow
in $\B_\Phi$ is an isomorphism if and only if the forgetful functor
to $\B$ sends it into an isomorphism, \eqref{proverab} follows.
Exactness of the forgetful functor is clear from the explicit description
of $\Ker$, $\Coker$.
\epr

For every $X\in \B$ the object $\Phi(X)$ is naturally equipped with a structure
of a coalgebra over $(\Phi,\eta,\mu)$; thus we get a functor $\Phitil:
\B \to \B_\Phi$ (the ``cofree coalgebra'' construction).
This functor is right adjoint to the forgetful functor $\B_\Phi\to \B$,
i.e. we have a natural isomorphism
\begin{equation}\label{alw}
\Hom_{\B_\Phi}(X,\Phitil(Y))\cong \Hom_\B(X,Y),
\end{equation}
where we omit the forgetful functor $\B_\Phi\to \B$ from our notation
(see \cite{ML} \S VI.2, Theorem 1).

\ssec{adjo}{Adjoint pairs and triality Theorem} Let $F^*:\A\to \B$
be a functor, and $F_*:\B\to \A$ be a right adjoint functor.
Then the composition $\Phi=F^*\circ F_*:\calB\to\calB$ is equipped
with natural transformations $\eta:\Phi\to \Id$ and $\mu:\Phi\to
\Phi\circ \Phi$, which together form a comonad.

The functor $F^*$ factors canonically  (cf. \cite{ML}, VI.3, Theorem 1)
through a functor \begin{equation}\label{canfu} \can:\calA\to \calB_\Phi.
\end{equation}

The next statement is a particular case of the general
``triality Theorem'' of Barr and Beck
(see \cite{ML}, \S VI.7, Theorem 1 (and Exercise 6))
which gives an explicit criterion for $\can$ to be an equivalence;
we include the proof since it is a bit shorter in our particular case.

%\begin{Thm}\label
\th{equiv}
Let $\A$, $\B$ be abelian categories, $F^*:\A\to \B$ be an
additive functor, and $F_*:\B\to \A$ be a right adjoint functor.
Assume that $F^*$ is exact and faithful. Then \eqref{canfu}
provides an equivalence $\can: \A\cong \B_\Phi$.
%\end{Thm}
\eth

\prf $\can$ is faithful and exact, because such is its composition
with the faithful exact forgetful functor $\B_\Phi\to \B$.

Let us check that $\can$ is a full imbedding, i.e.
\begin{equation}\label{fullim}
\Hom_\A(M,N) \iso \Hom_{\B_\Phi}(\can(M), \can(N)).
\end{equation}

First we claim that
\eqref{fullim} holds if $N=F_*(Y)$; moreover, this is true for any pair
of adjoint functors $F_*$, $F^*$ (not necessarily between
additive categories). Indeed, using \eqref{alw} we get
$$
\Hom_{\B_\Phi}(\can(M), \can(N))\cong
\Hom_{\B_\Phi}(\can(M),\Phitil(Y))\cong
\Hom_{\B}(F^*M,Y)\cong \Hom_\A(M,N),
$$
and it is immediate to see that the resulting isomorphism coincides
with the map induced by $\can$.

Thus to check \eqref{fullim} it suffices to see that any $N\in \A$
is a kernel of an arrow $F_*(X)\to F_*(Y)$ for some $X,Y\in \B$.
It is enough to find an injection $N\imbed F_*(X)$ (then apply the same
construction to its cokernel); but $F^*$ being faithful  implies
that the adjunction arrow $N\to F_*F^*(N)$ is an injection.

To see that $\can$ is  surjective on isomorphism classes of objects
it sufficies to prove that any object $X\in \B_\Phi$ is a subobject in
$\can (M)$ for some $M\in \A$ (then it is also a kernel of a morphism
$\phi:\can(M)\to \can(N) $ for some $M,N\in \A$; since we know already
that $\phi=\can(\phi')$ for some $\phi'$ we conclude that $X\cong
\can (\Ker(\phi'))$ by exactness of $\can$).
Now for $X\in \B_\Phi$ consider the adjunction arrow
$X\to \Phitil(X) = \can F_*(X)$
(coming from \eqref{alw}). We claim it is injective; indeed, since
the forgetful functor is exact and faithful, it is enough to see
that the corresponding arrow $h:
X\to \Phi(X)$ in $\B$ is injective;
it is in fact a split injection, because $\eta_X\circ h=id$ by the
definition of a coalgebra over a comonad.
\epr

\ssec{locdef}{Gluing terminology}The following definitions appear
(in a somewhat less economical notation) in \cite{kl} and~\cite{Po}.
(In \cite{Po}  left adjoint functors are used instead of right adjoint
ones; in~\cite{kl} it is assumed that both left and right adjoint
functors exist.)

Let $W$ be a finite set and let
 $\{\calB_w\}_{w\in W}$ be a collection of abelian categories,
and let $\B=\oplus_w \B_w$ be their product. For a functor
$\Phi:\B\to \B$ we will write $\Phi=(\Phi_{w_1, w_2})$, where
$\Phi_{w_1,w_2}:\B_{w_2}\to \B_{w_1}$.

\defe{gldef}
A {\it right gluing data} for $(\B_w)$ is a comonad
$(\Phi, \eta, \mu)$ on $\B=\oplus_w \B_w$ such that

1) $\Phi$ is (additive and) left-exact;

2) for each $w\in W$ the morphism $\Phi_{w,w}\to \Id_{\B_w}$
induced by $\eta$ is an isomorphism. \edefe

For a gluing data $(\B_w, \Phi,\eta,\mu)$ we call the category
$\B_\Phi$ the category {\it glued from $B_w$}.

\defe{locdefe} A right localization data for $(\B_w)$ is an abelian
category $\A$ together with a collection of exact functors
$F_w^*:\A\to \B_w$, such that $F_w^*$ has a right adjoint $F_{w*}$,
and the adjunction arrow $F_w^*\circ F_{w*}\to Id_{\B_w}$ is an
isomorphism.

A localization data is called faithful if the sum $\oplus_w
F_w^*:\A\to \oplusl_w \B_w$ is faithful. \edefe

\refl{24}, \reft{equiv} imply the following

%\begin{Cor}\label{glueeqv}
\cor{glueeqv}
If $(\A,\B_w, F_w^*)$ is a localization data, then $\Phi=(\Phi_{w_1,w_2})$,
$\Phi_{w_1,w_2}=F_{w_1*}F_{w_2}^*$ is a gluing data for $\B_w$.

We have a canonical exact functor from $\A$ to the glued category $\B_\Phi$;
it is an equivalence if and only if the localization data is faithful.
$\square$
%\end{Cor}
\ecor

Conversely, if $(\B_w, \Phi)$ is a gluing data, then the glued category
$\B_\Phi$ with the forgetful functor $\B_\Phi\to \oplus_w \B_w$
is a faithful localization data for $\B_w$.

\rem{} Let $F^*:\A\to \B$, $F_*:\B\to\A$ be a pair of adjoint
functors between abelian categories $\A$, $\B$, where $F^*$ is
exact. It is not difficult to show that  the adjunction arrow
$F^*\circ F_*\to \Id_\B$ is an isomorphism if and only if $F^*$
induces an equivalence $\A/\calC\cong \B$, where $\calC\subset \A$
is a Serre subcategory, and $\A/\calC$ is the Serre quotient.

Thus localization data amounts to the data of an abelian category
$\A$, and a  finite collection of Serre subcategories $\calC_w\subset
\A$, such that the projection $\A\to \A/\calC_w$ admits a right
adjoint. The localization data is faithful if the intersection
$\cap \calC_w$ is the zero category
(recall that a Serre subcategory is strictly full by definition).
Since any gluing data admits a unique (up to equivalence)
faithful localization data,
%%(namely, the glued category),
we see that a gluing data amounts to the data of an abelian category
$\A$, and a finite collection of Serre subcategories $\calC_w$ with
zero intersection, and such that the projection $\A\to \A/\calC_w$
has a right adjoint.

Assume that an abelian category $\A$ admits arbitrary direct sums,
and the functors of direct sums in $\A$ are exact.  Assume also
that $\A$ admits a set of generators.  Then it follows from
\cite{faith}, 5.51--5.53 that  for a Serre subcategory $\calC\subset \A$
the projection $\A\to \A/\calC$ admits a right adjoint if and only if
$\calC$ is closed under arbitrary direct sums.
\erem

%Let us also recall
%that if an abelian category $\A$ admits arbitrary
%direct sums, then
%
% (see e.g. \cite{faith}5.51-5.53

%----------------------------------------------------------------------------
\ssec{}{Examples}
%----------------------------------------------------------------------------
\subsubsection{}\label{sh} Let $X$ be a topological space and let
$\{U_w\}_{w\in W}$ be a collection of open subsets of $X$. Set
$\calA$ to be the category of sheaves of abelian groups on $X$,
$\calB_w$ be the category of sheaves of abelian groups on $U_w$; let
$F_w^*$ be the restriction functor and $F_*^w$ be the functor of
direct image of sheaves (this example explains our terminology).
%-----------------------------------------------------------------------

%----------------------------------------------------------------------------
\subsubsection{}\label{difex}
Let $X$ be a nonsingular quasi-affine complex algebraic
variety. Denote
by $D_X$ the sheaf of differential operators on $X$ and set
$\calD$ to be the algebra of global sections of $D_X$. Let
$W$ be any finite set of automorphisms of $\calD$ and set $\calA$ to be
the category of $D_X$-modules. Then every $w\in W$
defines a functor $M\mapsto M^w$ (twisting of the action by~$w$).
Set
\begin{itemize}
%\noindent $\bullet$
\item $\calB_w$ to be the category of quasi-coherent sheaves
of $D_X$-modules for any $w\in W$
%
%\noindent $\bullet$
\item $F^*_w(M)=D_X\underset{\calD}\ten M^w$
for every $M\in \calA$
%
%\noindent $\bullet$
\item $F_*^w(\calM)=\Gam(\calM)^{w^{-1}}$ for any
$w\in W$  and $\calM\in \calB_w$ (here $\Gam$ denotes the functor
of global sections for $D_X$-modules).
\end{itemize}
This is another example of a localization data whose particular case
is considered in  \refs{diff} below.
%---------------------------------------------------------------------------

\subsubsection{}
Let $A$ be a ring and $\A$ be the category of left $A$-modules.
For any idempotent element $e\in A$ consider the full subcategory
$\calC_e\subset\A$ whose objects are all $A$-modules $M$ such that
$eM=0$.  It is easy to see that $\calC_e$ is a Serre subcategory.

Consider the subring $eAe$ of the ring $A$.  It does not contain
the unit element of~$A$, but it has its own unit~$e$.  We claim
that the quotient category $\A/\calC_e$ can be identified with
the category of left $eAe$-modules, and the projection functor
sends an $A$-module $M$ to the $eAe$-module $eM$.
Indeed, let us set $B=eAe$ and denote by $\B$ the category of
left $B$-modules.  Let $F^*:\A\to\B$ be the functor $M\mapsto eM$.
We have
 $$
  F^*(M) \cong eA\otimes_A M \cong \Hom_A(Ae, M),
 $$
where $eA$ is considered as a $B$-$A$-bimodule and $Ae$ is
an $A$-$B$-bimodule.  Therefore, both left and right adjoint
functors to~$F^*$ exist and they can be computed as
 $$
  F_!(N) = Ae\otimes_{eAe} N, \quad
  F_*(N) = \Hom_{eAe}(eA, N).
 $$
It is easy to check that $F^*\circ F_*\cong\Id_\B\cong F^*\circ F_!$.

Now suppose that we are given a finite set of idempotents $e_w\in A$.
Let $\B_w$ denote the categories of left modules over the rings
$B_w=e_wAe_w$ and $F_w^*$, $F_{w*}$ be the above-defined functors.
This is obviously a localization data.  It is faithful
whenever $\sum_w Ae_wA=A$.

\subsubsection{}
 The following examples show how badly can the homological
dimension behave with respect to our gluing.  Let $A$ be the
associative algebra (over a field~$k$) generated by elements
$e_1$, $e_2$, $x_{12}$, and $x_{21}$ with the following relations:
$e_i$ are orthogonal idempotents, $e_1+e_2=1$, and
$e_ix_{jk}=\delta_{ij}x_{jk}$.
 Consider the quotient algebras $\tilde A=A/(x_{12}x_{21})$ and
$\hat A=\tilde A/(x_{21}x_{12})$.
 The last two algebras are finite-dimensional: $\dim_k \tilde A=5$
and $\dim_k \hat A=4$.

It is easy to see that the algebra $\hat A$ has infinite homological
dimension.  On the other hand, we have $e_1\hat A e_1=k$ and
$e_2\hat A e_2=k$.  So an abelian category of infinite homological
dimension is glued out of two copies of the category of vector spaces.

The algebra $\tilde A$ has homological dimension~2.  On the other hand,
we have $e_1\tilde Ae_1=k$ and $e_2\tilde Ae_2 \simeq k[y]/y^2$, where
$y=x_{21}x_{12}$.  Thus a category of homological dimension~2 is glued
out of the category of vector spaces and a category of infinite
homological dimension.
%%of one category of infinite homological dimension and one copy
%%of the category of vector spaces.

\sec{diff}{Gluing of $D$-modules on the basic affine space}
%----------------------------------------------------------------------------
\ssec{not}{Notation}In this section we will deal with
algebraic varieties over
$\CC$. For any such variety $Y$, the symbol $\calD(Y)$ will denote
the algebra of global differential operators on $Y$.
%---------------------------------------------------------------------------
\ssec{notgroup}{}
Let $G$ be a semisimple
simply-connected algebraic group over $\CC$;
let $U\subset G$ be a maximal unipotent subgroup, and let $X$ denote
the homogeneous space $G/U$. Then $U$ is the unipotent radical of a
Borel subgroup
$B\subset G$, and $X$ is known as the basic affine space of $G$.
The variety $X$ is a quasi-affine. Let $\calD=\calD(X)$ denote
the algebra of global
differential operators on $X$, and let $\calO$ denote the algebra of regular
functions on $X$.

For a simple root $\alpha$ of $G$
let $P_\alpha\subset G$ be the minimal parabolic subgroup
of type $\alpha$ containing $B$. Let $B_\alpha=[P_\alp,P_\alp]$
be the commutator
subgroup of $P_\alpha$,
and set $X_\alpha:=G/B_\alpha$. We have the obvious projection
of homogeneous spaces $\pi _\alpha:
X\rightarrow X_\alpha$.
It is a fibration with the fiber $B_\alpha /U=\AA^2-\{0\}$ (here
$\AA^2$ denotes the affine plane).

Let $\overline \pi _\alpha :\overline X^\alpha \to X_\alpha$ be the relative
affine completion of the morphism $\pi _\alpha$. (So  $\overline \pi _\alpha$
is the affine morphism corresponding to the sheaf of algebras $\pi_{\alpha *}
(\calO_X)$ on $X_\alpha$.) Then  $\overline \pi _\alpha$ has the structure
of a 2-dimensional vector bundle; $X$ is identified with the complement to
the zero-section in $\overline X^\alpha$. The $G$-action on $X$ obviously
extends to $\overline X^\alpha$; moreover, it is easy to see that the
determinant of the vector bundle $\overline \pi _\alpha$ admits a unique
(up to a constant) $G$-invariant trivialization, i.e.
$\overline \pi _\alpha$ admits unique up to a constant $G$-invariant
fiberwise symplectic form $\ome_\alp$ (cf. \cite{kl} for more
details).
In what follows we will
fix these forms for every simple root $\alp$.

Recall that for any symplectic vector bundle $p:E\to Z$
we have a canonical automorphism $\bfF_p$ of
the sheaf of algebras $p_*\calD_E$ and, in particular, of the ring
$\calD_E$ of global differential operators on $E$,
called  Fourier transform.
In particular, we get a canonical automorphism $F_\alpha=
F_{\overline \pi_\alpha}$ of the ring
$\calD(\overline X^\alpha)=\calD(X)
=\calD$
(the first equality follows from the fact that $\overline X^\alpha -X$
has codimension 2 in $\overline X^\alpha$).

\prop{action}(cf. \cite{GG})
Let $W$ be the Weyl group, $s_\alpha \in W$ be the simple
reflection of type $\alpha$. Then the assignment $s_\alpha \mapsto
\bfF_\alpha$
extends to a homomorphism $w\mapsto \bfF_w$ of the Weyl group to the group
of automorphisms of $\calD$.
\eprop
\prf Let $G_\RR$ denote the group $\RR$-points of a split real form of $G$.
Let also $U_\RR\subset G_\RR$ be the group of points of a maximal unipotent
subgroup of $G$ defined over $\RR$. Then the manifold $X_\RR=G_\RR/U_\RR$ admits a unique
up to a constant $G_\RR$-invariant measure which has unique smooth extension to
every $\oX^\alp_\RR$. Let $L^2(X_\RR)$ denote the
space of $L^2$-functions on $X_\RR$ with respect to this measure.
Since for every $\alp$ as above we have $L^2(X_\RR)=L^2(\oX^\alp_\RR)$ it
follows that we have well defined unitary operators $F_\alp$ acting on
$L^2(X_\RR)$ (Fourier transform along the fibers of $\opi_\alp$).
We claim now that these operators $F_\alp$ define an action of $W$
on $L^2(X_\RR)$. Indeed, this is proved in \cite{Ka} when $\RR$ is replaced
by a non-archimedian local field and in the archimedian case is essentially
a word-by-word repetition. For every $w\in W$ we denote by $F_w$ the
corresponding unitary automorphism of $L^2(X_\RR)$. The operators $F_w$
commute with the natural action of $G_\RR$ on $L^2(X_\RR)$.

Set now
\eq{}
\calS(X_\RR)=\{ f\in C^{\infty}(X_\RR)|\ d(f)\in L^2(X_\RR)\forall
d\in\calD\}
\end{equation}
(here by $C^\infty(X_\RR)$ we mean the space of {\it complex valued}
$C^\infty$-functions).
It is clear that $\calS(X_\RR)$ is a dense subspace of $L^2(X_\RR)$ (since
it contains the dense subspace of $C^\infty$-functions with compact
support).

Let us show that $\calS(X_\RR)$ is invariant with respect to the
operators $F_w$. Indeed, every $f\in\calS(X_\RR)$ is a
$C^\infty$-vector in the $G_\RR$ representation $L^2(X_\RR)$.
Therefore, since every $F_w$ commutes with $G_\RR$ it follows that
$F_w(f)$ is again a $C^\infty$-vector with respect to $G_\RR$.
Hence for every $d\in\calD$ the function $d(F_w(f))$ makes sense.
Moreover, it is easy to see that for every
$d\in\calD,f\in\calS(X_\RR)$ and for every simple root $\alp$ of
$G$ we have \eq{sopryazhenie}
F_\alp(\bfF_\alp(d)(f))=d(F_\alp(f)).
\end{equation}
Hence $d(F_\alp(f))\in L^2(X_\RR)$ which implies that
$F_\alp(f)\in\calS(X_\RR)$. Hence $\calS(X_\RR)$ is invariant with respect
to $F_\alp$'s and therefore it is also invariant with respect to all $F_w$.

It is clear that $\S(X_\RR)$ is a faithful module over $\calD$. Moreover
\refe{sopryazhenie} implies that in the space $\End \calS(X_\RR)$ we have
the equality
\eq{}
\bfF_\alp(d)=F_\alp\circ d\circ F_\alp= F_\alp^{-1}\circ d \circ F_\alp
\end{equation}
Clearly, this implies our claim.
\epr

\medskip
\rem{} One can also give an algebraic proof of
\refp{action} (the braid relations can be verified
using the analogue of the Radon transform associated
to any $w\in W$). We do not present details in this paper.
\erem
%-----------------------------------------------------------------------------
%-------------------------------
\ssec{}{}Let $H=B/U$ be the Cartan group of $G$ and let
$\grh=\text{Lie}(H)$.
Then $X$ carries a natural action of $H$,
commuting with the $G$-action (it comes from the action
of $B$ on $G$ by right translations). So $H$ acts in a locally finite
way on the rings $\calO, \calD$,
i.e. these rings are $\Lambda$-graded, where $\Lambda=\Hom (H, \CC^*)$ is the
weight lattice of $G$. Let $\calO^\lambda$ (resp. $\calD^\lam$)
denote the graded component of $\calO$ (resp. of $\calD$) of degree
$\lambda$. Let also $\Lambda ^+\subset \Lambda$ be the set of dominant
weights. We denote by $\rho\in \Lam^+$ the half-sum of the positive
roots.

Note that every element $h\in \grh$ defines a $G$-invariant vector
field on $X$. This defines an embedding of algebras
$U(\grh)\hookrightarrow \calD$, where $U(\grh)$ is the universal
enveloping algebra of $\grh$.

The following Lemma is an immediate consequence of the definitions
and \refp{action}.
\lem{fact}
\begin{enumerate}
\item The operators $\bfF_w$ commute with the $G$-action on $\calD$.
\item  For every  $h\in \grh$ we have
\eq{}
\bfF_w (h)=w(h)+\langle w(h)-h ,\rho \rangle 1
\end{equation}
(Here $\langle\;,\;\rangle$ is the natural pairing between $\grh$ and
$\Lambda$).
\end{enumerate}
In particular we have $\bfF_w(\calD(\lambda)^\mu)=\calD(\lambda)^{w(\mu)}$
\elem
%--------------------------------------------------------------------------
\ssec{main-prem}{A non-vanishing statement}Here we formulate
our main computational result (it is proved in \refss{pruf} below).
Let us denote by $\M$ be the category of modules over the ring $\dg$;
and let as before $\calM(D_X)$ be the category of $D$-modules on $X$.
We have the functor
of global sections $\Gamma :\calM(D_X)\to \M$, and the left-adjoint
functor $L:\M\to \calM(D_X)$.
For an affine open $U\subset X$ and $M\in \M$
we have:
\eq{}
 L(M)(U)=M\underset{\calD}\otimes \calD(U)=
M\underset{\calO}\otimes \calO (U)
\end{equation}
It is easy to see that $L$ is an exact functor, and that $L\circ \Gam\iso Id$.

For any $w\in W$ we will use the same notation $\bfF_w$ for an automorphism
of an associative ring $\dg$
and the corresponding auto-equivalence of the category
$\M$ of modules over $\calD$. Set
%--------------------------------------------------------------
\eq{}
L_w := L\circ \bfF_w,\quad\Gam_w:=\bfF_{w^{-1}}\circ \Gam
\end{equation}
%-------------------------------------------------------------------
%Then we see that the category $\M$ together with $|W|$ copies
%of the category $\M(D_X)$ and the functors $(L_w, \Gam_w)_{w\in W}$
%form a localization data in the sense of \refss{locdef}.
%---------------------------------------------------------------------------
\th{mainr}
For any $M\in Ob\M$, $M\not = 0$ there exists $w\in W$
such that $L_w (M)\not = 0$.
\eth
%-------------------------------------------------------------------------
\ssec{diff}{The structure of $\calM$}The statements formulated in the
Introduction are immediate consequences of \reft{mainr}

{\it Proof of  \reft{one}.}
If $I_1\subset I_2\subset \cdots $ is a chain of left ideals
in $\calD$ then $L_w(I_n)$ is a chain of
sub $D$-modules in the free $D$-module $D_X$. Since the category
of $D$-modules on a smooth variety is Noetherian, this chain stabilizes
for all $w\in W$; thus $L_w(I_{n+1}/I_n)=0$ for large $n$, hence $I_{n+1}
=I_n$ by \reft{mainr}. We checked that $\calD$ is left Noetherian; since
$\calD^{\mathrm {op}}\cong \calD$, because $X$ has a non-vanishing volume form,
we see that  $\calD$ is also right Noetherian. $\square$

\medskip

{\it Proof of \reft{two}.}
 It is enough  to show that $\Ext^i(M,\dg)=0$ for $i>
\dim X$ and any finitely generated $M\in  \calM$; and that
$\Ext^{\dim X}(M,\dg)\ne 0$ for some $M\in \calM$. Now,
$\Ext^i(M,\dg)$ may be considered as a right $\dg$-module by means
of the right action of $\dg$ on itself. Therefore, it is enough to
show that for any $w\in W$ one has $L_w(\Ext^i(M,\dg))=0$
(\reft{mainr} is valid, of course, also for right $\dg$-modules).
Suppose that this is not so, i.e. that there exists $w\in W$ such
that $L_w(\Ext^i(M,\dg))\neq 0$. We may assume, without loss of
generality, that $w=1$. For a finitely generated projective object
$P\in \calM$ we have a canonical isomorphism $L(\Hom(P,\dg))\cong
\calH om (L(P), D_X)$, where $\calH om$ is the sheaf of Hom's;
hence for any finitely generated $M\in \calM$ we have
$L(\Ext^i(M,\dg))\simeq \calE xt^i(L(M),D_X)$ where $\calE
xt^i(L(M),D_X)$ is the corresponding sheaf of Ext's. This means
that for any affine open subset $U$ of $X$ one has \eq{}
\Gam(U,L(\Ext^i(M,\calD)))=\Ext^i_{\calD(U)}(M,\calD(U)).
\end{equation}
But the  right hand side of this equality vanishes when $i>\dim
X=\dim U$, since for a non-singular affine variety $U$ the algebra
$\calD(U)$ has homological dimension equal to $\dim U$. Also, for
$i=\dim X$ it is non-zero, for example if $M=\Gamma(\calO)$, so
that $L(M)\cong \calO$.

\th{three}  $|W|$ copies of the
category $\M(D_X)$ (indexed by $W$)
together with  functors  $\Phi_{w_1,w_2}=\Gamma_{w_1}^*\circ L_{w_2*}$
and natural tranformations $\Phi_{w_1,w_3}\to \Phi_{w_1,w_2}\circ \Phi_{w_2,
w_3}$, $\Phi_{w,w}\to Id$ arising from adjointness of $L_w$ and $\Gamma_w$
form a gluing data.
The glued category is naturally equivalent to $\M$.
\eth

\prf \reft{mainr} shows that categories
 $\A=\calM$, \;  $\B_w = \M(D_X)$ for $w\in
W$; and functors $F_w^*=L_w:\A\to \B_w$; \;  $F_{w*}=\Gamma_w:\B_w\to \A$
form  a faithful localization data.
Thus the statement follows from \refc{glueeqv}. \epr

%----------------------------------------------------------------------------

 \ssec{pruf}{Proof of  \reft{mainr}} Let $\overline X=\Spec
\calO$ be the affine completion of $X$. We have the obvious open
embedding $X\hookrightarrow \oX$. Let $I\subset \calO$ be the
ideal of functions vanishing on $\oX\backslash X$. Then for
$M\in \M$
we have $L(M)=0$ if and only if  $I$ acts on $M$ locally
nilpotently.

% Let us recall a few well-known facts on the structure of the involved
%objects.

The following result is well-known.
%--------------------------------------------------------------------
\lem{borel-weil}
\begin{enumerate}
\item (Bott-Borel-Weil) $\calO^\lambda$ is  the irreducible $G$-module
of highest weight $\lambda$ when $\lambda\in \Lam^+ $ and
$\calO^\lambda = 0$ otherwise.
%-----------------------------------------------------------------
\item  The ideal $I$ is generated by $\calO^\rho\subset \calO$.
\end{enumerate}
\elem
%--------------------------------------------------------------------
We now proceed to the proof of \reft{mainr}. Assume that
$M\in Ob\M$ is such that $L(\bfF_w(M))=0$ for all $w\in W$. Then $\bfF_w(I)$
acts on $M$ locally nilpotently for all $w$. Fix $m\in M$, $m\not = 0$.
Then from the second statement of \refl{borel-weil}
we see that for some $n$ we have
$\bfF_w((\calO^\rho)^n)m=
\bfF_w (\calO^{n\rho})m=0$ for all $w$. So \reft{mainr} follows from
the following result. Let $w_0\in W$ denote the longest element. For any
$\lam\in\Lam$ we set $\lam^\vee=w_0(\lam)$.
%--------------------------------------------------------------------------
\prop{one} For any dominant weight $\lambda$ the left ideal
in $\dg$ generated by the space
$\sum _{w\in W} \bfF_w(A^{\lambda ^\vee {\;}} )$ contains 1.
\eprop
%--------------------------------------------------------------------------
For a $G$-module $V$, let $V(\mu)$ denote the isotypic part of $V$
corresponding to the irreducible $G$-module of highest weight $\mu$.
We start the proof of \refp{one} with the following
%----------------------------------------------------------------------------
\lem{dzero}
(cf. \cite{sh} for a different proof). We have: $\calD(0)^0 =  U(\grh)$.
\elem

\prf
It is obvious that $U(\grh) \subset \calD(0)^0$.
Let us prove the inverse inclusion.

Let $d\in \calD(0)^0$.
Since $\calO^\lambda$ is irreducible for any $\lambda$, we have
$d|_{\calO^\lambda}=c_\lambda$ for some $c_\lam \in \CC$.

Let us prove that the function $c:\Lambda^+ \to \CC$, $c(\lam):=
c_\lam$ is polynomial.
For any $\lambda \in \Lambda ^+$ consider the operator $T_\lambda$
acting on the functions on $\Lambda ^+$ which is defined by:
$T_\lambda (f)(\mu):=f(\lambda +\mu)-f(\mu)$.
Then for any set of functions $f_0,...,f_{n+1}$, where
$f_i\in \calO^{\lambda _i}$,
we have:
%----------------------------------------------------------------------
\eq{}
[f_0,[f_1,...[f_n, d]..]](f_{n+1})=T_{\lambda_0}\circ
T_{\lambda _1}\circ ...\circ T_{\lambda _n} (c) (\lambda _{n+1})
f_0 f_1 ... f_{n+1}
\end{equation}
%------------------------------------------------------------------------
In particular, if $d$ is a differential operator of order
$n$, then $T_{\lambda_0}\circ
T_{\lambda _1}\circ ...\circ T_{\lambda _n} (c)=0$ for any
$\lambda _0, \lambda _1,...,\lambda _n\in \Lambda ^+$. But it is well known
that the latter property implies that the function $c$ is polynomial.

It is clear that for every polynomial function $\phi :\Lambda ^+
\to \C$ there exists an element $u_\phi \in U(\h)$ such that
$u_\phi |_{\calO^\lambda}=\phi (\lambda)$. Let us take $u=u_c$.
Then $d (f)=u(f)$ for any $f \in \calO$. Hence $d=u \in U(\h)$
which finishes the proof.
\epr

The plan of the proof of \refp{one} is as follows.
Let $w_0\in W$ as before be the longest element.
Consider the map $m:\calO^\lambda \otimes \calO^{\lambda \galochka}\to
\calD$, defined by $m(f\otimes g):=f\cdot \bfF_{w_0}(g)$. Then $m$ lands to
$\calD^{\lambda +w_0(\lambda \galochka)}=\calD^0$.

Let $C_\lambda \in \calO^\lambda \otimes \calO^{\lambda \galochka}$
be the unique
(up to a constant) $G$-invariant element. (Recall that
$\calO^\lambda \otimes \calO^{\lambda \galochka}\simeq V(\lambda)\otimes
V(\lambda \galochka)\simeq V(\lambda)\otimes V(\lambda)^*$.)
Then $m(C_\lambda)\in (\calD^0)^G = \calD(0)^0=U(\h)$ by \refl{dzero}.

Let $P_\lambda$ denote $m(C_\lambda)$. It is clear that for any $w\in W$
the element $\bfF_w(P_\lambda)$
lies in the left ideal generated by
$\sum_{w\in W} \bfF_w(\calO^{\lambda ^\vee {\;}} )$.

On the other hand  we will prove the following
%----------------------------------------------------------------------------
\prop{two}The element $P_\lambda \in U(\grh)$ is of the form
\eq{}
0\neq P_{\lambda}=\text{const}\prod
(\alpha^{\vee} _i +\langle \alpha^{\vee} _i , \rho \rangle -l_i)
\end{equation}
where $\alpha^{\vee} _i \in \h$ are positive coroots, $l_i$ are positive
integers
and $\text{const}\in \CC^*$ is a non-zero constant.
\eprop
%----------------------------------------------------------------------
\refp{two} implies the following
%--------------------------------------------------------------------------
\cor{corollary}The ideal in $U(\grh)$ generated by $\bfF_w(P_\lambda)$
for all $w \in W$ contains 1.
\ecor
%----------------------------------------------------------------

\noindent
{\it Proof of the Corollary.} By Hilbert Nullstellensatz it suffices
to prove that for any point $x\in \h ^*$ there exists $w\in W$
such that $\bfF_w(P_\lambda)(x)\not = 0$. (Here we identified $U(\h)=S(\h)$
with the algebra of polynomial functions on $\h ^*$.)
It is enough to take $w$ such that $w^{-1}(x-\rho)$ is an
antidominant weight. $\square$

\refc{corollary} obviously implies \refp{one}.

Let us now prove \refp{two}. We will prove a more precise
%---------------------------------------------------------------------------
\prop{three}We have:
\eq{mainformula}
0\neq  P_{\lambda} = \text{const}\prod _{\alpha^{\vee} \in \Sigma ^+}
\prod _{i=1}^{\langle \alpha^{\vee},
\lambda \rangle} (\alpha^{\vee} + \langle \alpha^{\vee}, \rho \rangle -i)
\end{equation}
%------------------------------------------------------------
where $\Sigma ^+\subset \h$ is the set of positive coroots and
$\text{const}\in \CC^*$.
\eprop
%--------------------------------------------------------------------------
The first step towards the proof is the following
\lem{order}We have
\eq{}
\order(P_\lambda) \leq \langle \lambda , 2 \rho \rangle
\end{equation}
Here order denotes the standard order of a differential operator.
\elem
\prf
Let $p:E\to Y$ as before be a symplectic vector bundle.
Let $d\in \calD(E)$. Suppose that $d$ transforms by a character $t\mapsto
t^n$ under the natural action of the group $\CC^*$ on $\calD(E)$
(coming
from the action of $\CC^*$ on $E$ by dilatations). Then it is easy to
check that
\eq{}
\order(\bfF_p(d))=\order(d)+n
\end{equation}
For $d\in \calD^\lambda$ we know that $\bfF_w(d) \in
\calD^{w(\lambda)}$ .
So dilatations in the fibers of the vector bundle
$\overline \pi _\alpha$
act on $\bfF_w (d)$ by the character
$t\mapsto t^{\langle \alpha^{\vee} , w(\lambda)\rangle}
= t^{\langle w(\lambda) - s_\alpha w (\lambda), \rho \rangle}$.
By induction on the length of $w$ we deduce that $\order(\bfF_w(d))=
\order(d)+\langle \lambda - w(\lambda), \rho \rangle$. In particular,
for any $\mu \in \Lambda ^+$ and $f \in \calO^\mu$ we have:
$\order(F_{w_0}(f))=
\langle \mu - w_0(\mu), \rho \rangle=\langle \mu ,2 \rho \rangle$.
The lemma follows.
\epr
%---------------------------------------------------------------------------
\lem{divides}
For $\alpha^{\vee} \in \Sigma ^+$ and an integer $i$ such that
$0\leq i  < \langle \alpha^{\vee} , \lambda \rangle $ the element
$(\alpha^{\vee} + \langle \alpha^{\vee} , \rho \rangle -i -1)\in\grh$
divides $P_\lambda$.
\elem
\prf
Let us choose a simple coroot $r$ and $w\in W$ such that
$w(\alpha^{\vee})=r$. By \refl{fact} the statement of \refl{divides}
is equivalent
to saying that $(r-i)$ divides $\bfF_w(P_\lambda)$.

We are going to show that for $x \in \Lambda ^+$ such that
$ \langle r , x \rangle =i$ we have $\bfF_{w  }(P_\lambda)(x) =0$.
(Recall that we have identified $U(\h)$ with the algebra of polynomial
functions on $\h ^*$).
This is enough, since the set $\{x \in \Lambda ^+|
  \langle r , x \rangle =i\}$ is a Zariski dense subset of the hyperplane
$\{ x\in \grh^*|\ \langle r,x\rangle =i\}$.

The latter statement is equivalent to saying that $\bfF_w(P_\lambda)(f)=0$
for $f \in \calO^x$. Let us prove this.

We have $P_\lambda = \sum_i f_i \bfF_{w_0} (g_i)$ where $\{f_i\}$, $\{g_i\}$
are the dual bases of $\calO^\lambda$ and
$\calO^{\lambda \galochka}$ respectively.
So $\bfF_w(P_\lambda) = \sum_i \bfF_w(f_i) \bfF_{w w_0} (g_i)$.

We claim that $\bfF_{w w_0} (g_i) (f)=0$ for $f\in \calO^x$. Indeed,
$\bfF_{w w_0} (g_i) \in \calO^{w  (-\lambda)}$, hence
$\bfF_{w w_0}(g_i) (f) \in \calO^{x-w(\lambda)}=0$ since
$x-w(\lambda)$ is not dominant. (Recall that
$\langle r, x - w  (\lambda)\rangle = i - \langle \alpha^{\vee},
\lambda \rangle <0$ by the assumption).

This proves the lemma.
\epr
Now we are ready to prove \refp{three}. Let
RHS denote the right hand side of equality \refe{mainformula}. From
\refl{divides} it follows that RHS divides $P_\lambda$. From
\refl{order} we see that $\text{order}(P_\lambda)\leq \text{order(RHS)}$.

Since both $P_\lam$ and RHS are $G$-invariant and since
$\calO^G=\CC$ the equality will follow provided we know that $P_\lambda\neq 0$.

To check this take
the dual bases $\{f_i\}$, $\{g_i\}$ of $\calO^{\lam}$ and
$\calO^{\lam \galochka}$ respectively, compatible with the  weight
decomposition. Assume that $f_1$ is a highest weight vector, and
$g_1$ is a lowest weight vector. Then $\bfF_{w_0}(g_1)\neq 0$ which
implies that there exists $\mu\in \Lam^+$ and a highest weight vector
$\phi\in \calO^\mu$, such that
$\bfF_{w_0}g_1(\phi)\not = 0$  and
$\bfF_{w_0}(g_i)(\phi)=0$ for $i>1$. Hence  $P_\lambda (\phi)=
f_1 \bfF_{w_0}(g_1)(\phi)\not = 0$.

The proof is finished. $\square$

\end{document}